\documentclass[12pt]{amsart}
\input{epsf}
\usepackage{amssymb,latexsym,cite}
\topskip  \numberwithin{equation}{section}
\newtheorem{theorem}{Theorem}[section]

\newtheorem{remark}[theorem]{Remark}

\begin{document}

\pagenumbering{arabic}
\pagestyle{headings}

\newcommand{\DPB}[4]{P\beta_{#1}^{(#2)}(#3,#4)}

\title{Representation by degenerate Frobenius-Euler polynomials}
\author{Taekyun Kim}
\address{Department of Mathematics, Kwangwoon University, Seoul 139-701, Republic of Korea}
\email{tkkim@kw.ac.kr}

\author{Dae San Kim}
\address{Department of Mathematics, Sogang University, Seoul 121-742, Republic of Korea}
\email{dskim@sogang.ac.kr}


\begin{abstract}
The aim of this paper is to represent any polynomial in terms of the degenerate Frobenius-Euler polynomials and more generally of the higher-order degenerate Frobenius-Euler polynomials. We derive explicit formulas with the help of umbral calculus and illustrate our results with some examples.
\end{abstract}
\maketitle


\section{Introduction and preliminaries}

The aim of this paper is to derive formulas expressing any polynomial in terms of the degenerate Frobenius-Euler polynomials (see \eqref{9A}) with the help of umbral calculus (see Theorem 3.1) and to illustrate our results with some examples (see Section 5). This can be generalized to the higher-order degenerate Frobenius-Euler polynomials (see \eqref{10A}). Indeed, we deduce formulas for representing any polynomial in terms of the higher-order degenerate Frobenius-Euler polynomials again by using umbral calculus (see Theorem 4.1). As corollaries to these theorems, we obtain formulas for expressing any polynomial by the degenerate Euler (see Theorem 3.2) and the higher-order degenerate Euler polynomials (see Theorem 4.2). The contribution of this paper is the derivation of such formulas which have potential applications to finding many interesting polynomial identities, as illustrated in Section 5. \par

The following identity is a slight modification of  the formulas obtained in [12] which express any polynomial in terms of Bernoulli polynomials $B_{n}(x)$ defined by $\frac{t}{e^{t}-1} e^{xt}=\sum_{n=0}^{\infty}B_{n}(x)\frac{t^{n}}{n!}$:
\begin{align}
& \sum_{k=1}^{n-1}\frac{1}{2k\left(2n-2k\right)}B_{2k}\left(x\right)B_{2n-2k}\left(x\right)+\frac{2}{2n-1}B_{1}\left(x\right)B_{2n-1}\left(x\right)\label{1A}\\
= & \frac{1}{n}\sum_{k=1}^{n}\frac{1}{2k}\dbinom{2n}{2k}B_{2k}B_{2n-2k}\left(x\right)+\frac{1}{n}H_{2n-1}B_{2n}\left(x\right)\nonumber 
+\frac{2}{2n-1}B_{1}\left(x\right)B_{2n-1},\nonumber
\end{align}
where $n \ge 2$, and $H_{n}=1+\frac{1}{2}+ \cdots +\frac{1}{n}$. \par
Letting $x=0$ and $x=\frac{1}{2}$ in \eqref{1A} respectively give a slight variant of the Miki's identity and the Faber-Pandharipande-Zagier (FPZ) identity.

In 1998, Faber and Pandharipande found that the FPZ identity must be valid for certain conjectural relations between Hodge integrals in Gromov-Witten theory. In the appendix of [5], Zagier gave a proof of the FPZ identity.
Another proof of FPZ identity is given in Dunne-Schubert [4] by using the asymptotic expansion
of some special polynomials coming from the quantum field theory computations.
As to the Miki's identity, Miki [17] uses a formula for the Fermat quotient
$\frac{a^{p}-a}{p}$ modulo $p^{2}$, Shiratani-Yokoyama [21] utilizes $p$-adic analysis and Gessel [6] exploits two different expressions for Stirling numbers of the second kind $S_{2}\left(n,k\right)$. 
Here it should be stressed that these proofs of Miki's and FPZ identities are quite involved, while our proofs of Miki's and Faber-Pandharipande-Zagier identities follow from the simple formulas in [12], involving only derivatives and integrals of the given polynomials. \par
Many interesting identities have been derived by using these formulas and similar ones for Euler and Frobenius-Euler polynomials (see [8--13,15]). Some convolution identities for Frobenius-Euler polynomials are found in [7,19] by exploiting generating function methods and summation transform techniques. The list in the References are far from being exhaustive. However, the interested reader can easily find more related papers in the literature. Also, we should mention here that there are other ways of obtaining the same result as the one in \eqref{1A}. One of them is to use Fourier series expansion of the function obtained by extending by periodicity of period 1 of the polynomial function restricted to the interval $[0,1)$ (see [16]).

The outline of this paper is as follows. In Section 1, we recall some necessary facts that are needed throughout this paper. In Section 2, we go over umbral calculus briefly. In Section 3, we derive formulas expressing any polynomial in terms of the degenerate Frobenius-Euler polynomials. In Section 4, we derive formulas representing any polynomial in terms of the higher-order degenerate Frobenius-Euler polynomials. In Section 5, we illustrate our results with some examples. Finally, we conclude our paper in Section 6.

The Euler polynomials $E_n(x)$ are defined by
\begin{equation}\label{2A}
\frac{2}{e^t+1}e^{xt}=\sum_{n=0}^{\infty}E_{n}(x)\frac{t^n}{n!}.
\end{equation}
When $x=0$, $E_n=E_n(0)$ are called the Euler numbers. We observe that $E_n(x)=\sum_{j=0}^{n}\binom{n}{j}E_{n-j}x^j,\,\frac{d}{dx}E_n(x)=nE_{n-1}(x),\,E_{n}(x+1)+E_{n}(x)=2x^{n}$.
The first few terms of $E_n$ are given by:
\begin{align*}
&E_0=1,\,E_1=-\frac{1}{2},\,E_3=\frac{1}{4},\,E_5=-\frac{1}{2},\,E_7=\frac{17}{8},\,B_9=-\frac{31}{2},\dots;\\
&E_{2k}=0,\,\,(k \ge 1).
\end{align*}

More generally, for any nonnegative integer $r$, the Euler polynomials $E_n^{(r)}(x)$ of order $r$
are given by

\begin{equation}\label{3A}
\bigg(\frac{2}{e^t+1}\bigg)^{r}e^{xt}=\sum_{n=0}^{\infty}E_{n}^{(r)}(x)\frac{t^n}{n!}.
\end{equation}

The Frobenius-Euler polynomials $H_{n}(x|u)\,\,(u \ne 1)$ are defined by
\begin{equation}\label{4A}
\frac{1-u}{e^{t}-u}e^{xt}=\sum_{n=0}^{\infty}H_{n}(x|u)\frac{t^{n}}{n!}.
\end{equation}
When $x=0$, $H_{n}(u)=H_{n}(0|u)$ are called the Frobenius-Euler numbers. We observe that $H_n(x|u)=\sum_{j=0}^{n}\binom{n}{j}H_{n-j}(u)x^j,\,\frac{d}{dx}H_{n}(x|u)=nH_{n-1}(x|u),\,H_{n}(x+1|u)-uH_{n}(x|u) = (1-u)x^{n}$.
The first few terms of $H_n(u)$ are given by:
\begin{equation*}
H_{0}(u)=1,\, H_{1}(u)=-\frac{1}{1-u},\, H_{2}(u)= \frac{1+u}{(1-u)^{2}},\,  H_{3}(u)=-\frac{u^2+4u+1}{(1-u)^3},\dots.
\end{equation*}

More generally, for any nonnegative integer $r$, the Frobenius-Euler polynomials $H_n^{(r)}(x|u)\,\,(u \ne 1)$ of order $r$
are given by

\begin{equation}\label{5A}
\bigg(\frac{1-u}{e^t-u}\bigg)^{r}e^{xt}=\sum_{n=0}^{\infty}H_{n}^{(r)}(x|u)\frac{t^n}{n!}.
\end{equation}

For any nonzero real number $\lambda$, the degenerate exponentials are given by
\begin{align}
&e_{\lambda}^{x}(t)=(1+\lambda t)^{\frac{x}{\lambda}}=\sum_{n=0}^{\infty}(x)_{n,\lambda}\frac{t^n}{n!}, \label{6A}\\
&e_{\lambda}(t)=e_{\lambda}^{1}(t)=(1+\lambda t)^{\frac{1}{\lambda}}=\sum_{n=0}^{\infty}(1)_{n,\lambda}\frac{t^n}{n!}\nonumber.
\end{align}

Carlitz [1] introduced a degenerate version of the Euler polynomials $E_n(x)$, called the degenerate Euler polynomials and denoted by $\mathcal{E}_{n,\lambda}(x)$, which is given by
\begin{equation}\label{7A}
\frac{2}{e_{\lambda}(t)+1}e_{\lambda}^{x}(t)=\sum_{n=0}^{\infty}\mathcal{E}_{n,\lambda}(x)\frac{t^n}{n!}.
\end{equation}
For $x=0$, $\mathcal{E}_{n,\lambda}=\mathcal{E}_{n,\lambda}(0)$ are called the degenerate Euler numbers. \par
More generally, for any nonnegative integer $r$, the degenerate Euler polynomials $\mathcal{E}_{n,\lambda}^{(r)}(x)$ of order $r$
are given by
\begin{equation}\label{8A}
\bigg(\frac{2}{e_{\lambda}(t)+1}\bigg)^{r}e_{\lambda}^{x}(t)=\sum_{n=0}^{\infty}\mathcal{E}_{n,\lambda}^{(r)}(x)\frac{t^n}{n!}.
\end{equation}

A degenerate version of the Frobenius-Euler polynomials $H_{n}(x|u)\,\,(u \ne 1)$, called the degenerate Frobenius-Euler polynomials and denoted by $h_{n,\lambda}(x|u)$, is given by (see [14])
\begin{equation}\label{9A}
\frac{1-u}{e_{\lambda}(t)-u}e_{\lambda}^{x}(t)=\sum_{n=0}^{\infty}h_{n,\lambda}(x|u)\frac{t^n}{n!}.
\end{equation}
For $x=0$, $h_{n,\lambda}(u)=h_{n,\lambda}(0|u)$ are called the degenerate Frobenius-Euler numbers. Clearly, $h_{n,\lambda}(x|-1)= \mathcal{E}_{n,\lambda}(x)$.    \par
More generally, for any nonnegative integer $r$, the degenerate Frobenius-Euler polynomials $h_{n,\lambda}^{(r)}(x|u)\,\,(u \ne 1)$ of order $r$ are given by (see [14])
\begin{equation}\label{10A}
\bigg(\frac{1-u}{e_{\lambda}(t)-u}\bigg)^{r}e_{\lambda}^{x}(t)=\sum_{n=0}^{\infty}h_{n,\lambda}^{(r)}(x|u)\frac{t^n}{n!}.
\end{equation}

For $x=0$, $h_{n,\lambda}^{(r)}(u)=h_{n,\lambda}^{(r)}(0|u)$ are called the degenerate Frobenius-Euler numbers of order $r$. Obviously, $h_{n,\lambda}^{(r)}(x|-1)= \mathcal{E}_{n,\lambda}^{(r)}(x)$.    \par
We remark that $h_{n,\lambda}(x|u) \rightarrow H_{n}(x|u)$, and $h_{n,\lambda}^{(r)}(x|u) \rightarrow H_{n}^{(r)}(x|u)$, as $\lambda$ tends to $0$.

We recall some notations and facts about forward differences. Let $f$ be any complex-valued function of the real variable $x$. Then, for any real number $a$, the forward difference $\Delta_{a}$ is given by
\begin{equation}\label{11A}
\Delta_{a}f(x)=f(x+a)-f(x).
\end{equation}
If $a=1$, then we let
\begin{equation}\label{12A}
\Delta f(x)=\Delta_{1}f(x)=f(x+1)-f(x).
\end{equation}
We also need
\begin{equation*}
\tilde{\Delta}f(x)=f(x+1)+f(x).
\end{equation*}
It is necessary to note that
\begin{equation*}
\Delta_{a}\tilde{\Delta}=\tilde{\Delta}\Delta_{a}.
\end{equation*}

In general, the $n$th oder forward differences are given by
\begin{equation}\label{13A}
\Delta_{a}^{n}f(x)=\sum_{i=0}^{n}\binom{n}{i} (-1)^{n-i}f(x+ia).
\end{equation}
For $a=1$, we have
\begin{equation}\label{14A}
\Delta^{n}f(x)=\sum_{i=0}^{n}\binom{n}{i} (-1)^{n-i}f(x+i).
\end{equation}
It is easy to see that
\begin{equation}\label{15A}
\tilde{\Delta}^{n}f(x)=\sum_{i=0}^{n}\binom{n}{i}f(x+i).
\end{equation}
Finally, we recall that the Stirling numbers of the second kind $S_{2}(n,k)$ are given by
\begin{equation}\label{16A}
\frac{1}{k!}(e^{t}-1)^{k}=\sum_{n=k}^{\infty}S_{2}(n,k)\frac{t^{n}}{n!},\quad(k \ge 0).
\end{equation}

\section{Review of umbral calculus}
\vspace{0.5cm}
Here we will briefly go over very basic facts about umbral calculus. For more details on this, we recommend the reader to refer to [3,20,22].
Let $\mathbb{C}$ be the field of complex numbers. Then $\mathcal{F}$ denotes the algebra of formal power series in $t$ over $\mathbb{C}$, given by
\begin{displaymath}
 \mathcal{F}=\bigg\{f(t)=\sum_{k=0}^{\infty}a_{k}\frac{t^{k}}{k!}~\bigg|~a_{k}\in\mathbb{C}\bigg\},
\end{displaymath}
and $\mathbb{P}=\mathbb{C}[x]$ indicates the algebra of polynomials in $x$ with coefficients in $\mathbb{C}$. \par
The set of all linear functionals on $\mathbb{P}$ is a vector space as usual and denoted by $\mathbb{P}^{*}$. Let  $\langle L|p(x)\rangle$ denote the action of the linear functional $L$ on the polynomial $p(x)$. \par
For $f(t)\in\mathcal{F}$ with $\displaystyle f(t)=\sum_{k=0}^{\infty}a_{k}\frac{t^{k}}{k!}\displaystyle$, we define the linear functional on $\mathbb{P}$ by
\begin{equation}\label{1B}
\langle f(t)|x^{k}\rangle=a_{k}. 
\end{equation}
From \eqref{1B}, we note that
\begin{equation*}
 \langle t^{k}|x^{n}\rangle=n!\delta_{n,k},\quad(n,k\ge 0), 
\end{equation*}
where $\delta_{n,k}$ is the Kronecker's symbol. \par
Some remarkable linear functionals are as follows:
\begin{align}
&\langle e^{yt}|p(x) \rangle=p(y), \nonumber \\
&\langle e^{yt}-1|p(x) \rangle=p(y)-p(0), \label{2B} \\
& \bigg\langle \frac{e^{yt}-1}{t}\bigg |p(x) \bigg\rangle = \int_{0}^{y}p(u) du.\nonumber
\end{align}
Let
\begin{equation}\label{3B}
 f_{L}(t)=\sum_{k=0}^{\infty}\langle L|x^{k}\rangle\frac{t^{k}}{k!}.
\end{equation}
Then, by \eqref{1B} and \eqref{3B}, we get
\begin{displaymath}
    \langle f_{L}(t)|x^{n}\rangle=\langle L|x^{n}\rangle.
\end{displaymath}
That is, $f_{L}(t)=L$. Additionally, the map $L\longmapsto f_{L}(t)$ is a vector space isomorphism from $\mathbb{P}^{*}$ onto $\mathcal{F}$.\par  Henceforth, $\mathcal{F}$ denotes both the algebra of formal power series  in $t$ and the vector space of all linear functionals on $\mathbb{P}$. $\mathcal{F}$ is called the umbral algebra and the umbral calculus is the study of umbral algebra. 
For each nonnegative integer $k$, the differential operator $t^k$ on $\mathbb{P}$ is defined by
\begin{equation}\label{4B}
t^{k}x^n=\left\{\begin{array}{cc}
(n)_{k}x^{n-k}, & \textrm{if $k\le n$,}\\
0, & \textrm{if $k>n$.}
\end{array}\right. 
\end{equation}
Extending \eqref{4B} linearly, any power series
\begin{displaymath}
 f(t)=\sum_{k=0}^{\infty}\frac{a_{k}}{k!}t^{k}\in\mathcal{F}
\end{displaymath}
gives the differential operator on $\mathbb{P}$ defined by
\begin{equation}\label{5B}
 f(t)x^n=\sum_{k=0}^{n}\binom{n}{k}a_{k}x^{n-k},\quad(n\ge 0). 
\end{equation}
It should be observed that, for any formal power series $f(t)$ and any polynomial $p(x)$, we have
\begin{equation}\label{6B}
\langle f(t) | p(x) \rangle =\langle 1 | f(t)p(x) \rangle =f(t)p(x)|_{x=0}.
\end{equation}
Here we note that an element $f(t)$ of $\mathcal{F}$ is a formal power series, a linear functional and a differential  operator. Some notable differential operators are as follows: 
\begin{align}
&e^{yt}p(x)=p(x+y), \nonumber\\
&(e^{yt}-1)p(x)=p(x+y)-p(x), \label{7B}\\
&\frac{e^{yt}-1}{t}p(x)=\int_{x}^{x+y}p(u) du.\nonumber
\end{align}

The order $o(f(t))$ of the power series $f(t)(\ne 0)$ is the smallest integer for which $a_{k}$ does not vanish. If $o(f(t))=0$, then $f(t)$ is called an invertible series. If $o(f(t))=1$, then $f(t)$ is called a delta series. \par
For $f(t),g(t)\in\mathcal{F}$ with $o(f(t))=1$ and $o(g(t))=0$, there exists a unique sequence $s_{n}(x)$ (deg\,$s_{n}(x)=n$) of polynomials such that
\begin{equation} \label{8B}
\big\langle g(t)f(t)^{k}|s_{n}(x)\big\rangle=n!\delta_{n,k},\quad(n,k\ge 0).
\end{equation}
The sequence $s_{n}(x)$ is said to be the Sheffer sequence for $(g(t),f(t))$, which is denoted by $s_{n}(x)\sim (g(t),f(t))$. We observe from \eqref{8B} that 
\begin{equation}\label{9B}
s_{n}(x)=\frac{1}{g(t)}p_{n}(x),
\end{equation}
where $p_{n}(x)=g(t)s_{n}(x) \sim (1,f(t))$.\par
In particular, if $s_{n}(x) \sim (g(t),t)$, then $p_{n}(x)=x^n$, and hence 
\begin{equation}\label{10B}
s_{n}(x)=\frac{1}{g(t)}x^n.
\end{equation}

It is well known that $s_{n}(x)\sim (g(t),f(t))$ if and only if
\begin{equation}\label{11B}
\frac{1}{g\big(\overline{f}(t)\big)}e^{x\overline{f}(t)}=\sum_{k=0}^{\infty}\frac{s_{k}(x)}{k!}t^{k}, 
\end{equation}
for all $x\in\mathbb{C}$, where $\overline{f}(t)$ is the compositional inverse of $f(t)$ such that $\overline{f}(f(t))=f(\overline{f}(t))=t$. \par

The following equations \eqref{12B}, \eqref{13B}, and \eqref{14B} are equivalent to the fact that  $s_{n}\left(x\right)$ is Sheffer for $\left(g\left(t\right),f\left(t\right)\right)$, for some invertible $g(t)$: 
\begin{align}
f\left(t\right)s_{n}\left(x\right)&=ns_{n-1}\left(x\right),\quad\left(n\ge0\right),\label{12B}\\
s_{n}\left(x+y\right)&=\sum_{j=0}^{n}\binom{n}{j}s_{j}\left(x\right)p_{n-j}\left(y\right),\label{13B}
\end{align}
with $p_{n}\left(x\right)=g\left(t\right)s_{n}\left(x\right),$
\begin{equation}\label{14B}
s_{n}\left(x\right)=\sum_{j=0}^{n}\frac{1}{j!}\big\langle{g\left(\overline{f}\left(t\right)\right)^{-1}\overline{f}\left(t\right)^{j}}\big |{x^{n}\big\rangle}x^{j}.
\end{equation}

For $s_{n}(x)\sim(g(t),f(t))$, and $r_{n}(x)\sim(h(t),l(t))$, we have 
\begin{equation}\label{15B}
s_{n}\left(x\right)=\sum_{k=0}^{n}C_{n,k}r_{k}\left(x\right),\quad\left(n\ge0\right),
\end{equation}
where 
\begin{equation}\label{16B}
C_{n,k}=\frac{1}{k!}\bigg\langle{\frac{h\left(\overline{f}\left(t\right)\right)}{g\left(\overline{f}\left(t\right)\right)}l\left(\overline{f}\left(t\right)\right)^{k}}\Big |{x^{n}}\bigg \rangle.
\end{equation}

\section{Representation by degenerate Frobenius-Euler polynomials}

Our goal here is to find formulas expressing any polynomial in terms of the degenerate Frobenius-Euler polynomials $h_{n,\lambda}(x|u)$. 

From \eqref{9A}, \eqref{6A} and \eqref{11B}, we first observe that 
\begin{align}
h_{n,\lambda}(x|u) &\sim \big(g(t)=\frac{e^t -u}{1-u}, f(t)=\frac{1}{\lambda}(e^{\lambda t}-1)\big), \label{1C}\\
& (x)_{n,\lambda} \sim (1, f(t)=\frac{1}{\lambda}(e^{\lambda t}-1)).\label{2C}
\end{align}
From \eqref{11A}, \eqref{7B}, \eqref{12B}, \eqref{1C} and \eqref{2C}, we note that
\begin{align}
f(t)h_{n,\lambda}(x|u)&=nh_{n-1,\lambda}(x|u)=\frac{1}{\lambda}(e^{\lambda t}-1)h_{n,\lambda}(x|u)=\frac{1}{\lambda}\Delta_{\lambda}h_{n,\lambda}(x|u), \label{3C}\\
&f(t)(x)_{n,\lambda}=n(x)_{n-1,\lambda}.\label{4C}
\end{align}
It is immediate to see from \eqref{9A} that
\begin{align}
&h_{n,\lambda}(x+1|u)-uh_{n,\lambda}(x|u)=(1-u)(x)_{n,\lambda}, \label{5C}\\
&h_{n,\lambda}(1|u)-uh_{n,\lambda}(u)=(1-u)\delta_{n,0},\label{6C}
\end{align}
where $\delta_{n,0}$ is the Kronecker's delta. \par
Now, we assume that $p(x) \in \mathbb{C}[x]$ has degree $n$, and write $p(x)=\sum_{k=0}^{n}a_{k}h_{n,\lambda}(x|u)$.
For a fixed $u$, let $a(x)=p(x+1)-up(x)$. Then, from \eqref{5C}, we have
\begin{align}
a(x)&=\sum_{k=0}^{n}a_{k}(h_{k,\lambda}(x+1|u)-uh_{k,\lambda}(x|u)) \label{7C}\\
&=(1-u)\sum_{k=0}^{n}a_{k}(x)_{k,\lambda}\nonumber
\end{align}
For $r \ge 0$, from \eqref{7C} and \eqref{4C} we obtain
\begin{align}
(f(t))^{r}a(x)&=(f(t))^{r}(1-u)\sum_{k=0}^{n}a_k(x)_{k,\lambda} \label{8C} \\
&=(1-u)\sum_{k=r}^{n}k(k-1) \cdots (k-r+1)a_{k}(x)_{k-r,\lambda}.\nonumber
\end{align}
Letting $x=0$ in \eqref{8C}, we finally get
\begin{equation}\label{9C}
a_{r}=\frac{1}{(1-u) r!}(f(t))^{r}a(x)|_{x=0}=\frac{1}{(1-u)r!}\langle (f(t))^{r}|a(x) \rangle,\,\,(r \ge 0),
\end{equation}
An alternative expression of \eqref{9C} is given by\\
\begin{align}
a_r&=\frac{1}{(1-u)r! \lambda^{r}}\Delta_{\lambda}^{r}a(x)|_{x=0} \nonumber\\
&=\frac{1}{(1-u)r! \lambda^{r}}(\Delta_{\lambda}^{r}p(1)-u\Delta_{\lambda}^{r}p(0)) \label{10C}\\
&=\frac{1}{(1-u)r! \lambda^{r}}\big(\Delta_{\lambda}^{r}\tilde{\Delta} p(0)+(1-u)\Delta_{\lambda}^{r}p(0)\big), \nonumber
\end{align}
since $f(t)a(x)=\frac{1}{\lambda}(e^{\lambda t}-1)a(x)=\frac{1}{\lambda} \Delta_{\lambda}a(x)$. \par
From \eqref{13A}, we have another alternative expression of \eqref{9C} which is given by
\begin{align}
a_r&=\frac{1}{(1-u)r! \lambda^{r}}\Delta_{\lambda}^{r}a(x)|_{x=0} \nonumber\\
&=\frac{1}{(1-u)r! \lambda^{r}}\sum_{k=0}^{r}\binom{r}{k}(-1)^{r-k}a(x+k \lambda)|_{x=0} \label{11C}\\
&=\frac{1}{(1-u)r! \lambda^{r}}\sum_{k=0}^{r}\binom{r}{k}(-1)^{r-k}(p(1+k \lambda)-up(k \lambda)).\nonumber
\end{align}
From \eqref{16A}, we get yet another expression of $a_{r}$ as follows:
\begin{align}
a_{r}&=\frac{1}{(1-u) \lambda^{r}}\frac{1}{r!}(e^{\lambda t}-1)^{r}a(x)|_{x=0} \label{12C}\\
&=\frac{1}{(1-u)}\sum_{l=r}^{n}S_{2}(l,r)\frac{\lambda^{l-r}}{l!}\big(p^{(l)}(1)-up^{(l)}(0)\big),\nonumber
\end{align}
where $p^{(l)}(x)=(\frac{d}{dx})^{l}p(x)$. \par
Finally, from \eqref{9C}--\eqref{12C}, we get the following theorem.

\begin{theorem}
Let $p(x) \in \mathbb{C}[x], with\,\, \mathrm{deg}\, p(x)=n$, and let $f(t)=\frac{1}{\lambda}(e^{\lambda t}-1)$.  Then we have $p(x)=\sum_{k=0}^{n}a_k h_{k,\lambda}(x|u)$, where 
\begin{align*}
&a_{k}=\frac{1}{(1-u)k!}(f(t))^{k}(p(x+1)-up(x))|_{x=0}\\
&=\frac{1}{(1-u)k! \lambda^{k}}\big \langle \big(e^{\lambda t}-1\big)^{k} \big | p(x+1)-up(x)\big \rangle \\
&=\frac{1}{(1-u)k! \lambda^{k}}\big(\Delta_{\lambda}^{r}\tilde{\Delta} p(0)+(1-u)\Delta_{\lambda}^{r}p(0)\big) \\
&=\frac{1}{(1-u)k! \lambda^{k}}\sum_{j=0}^{k}\binom{k}{j}(-1)^{k-j}(p(1+j \lambda)-up(j \lambda)) \\
&=\frac{1}{(1-u)}\sum_{l=r}^{n}S_{2}(l,r)\frac{\lambda^{l-r}}{l!}\big(p^{(l)}(1)-up^{(l)}(0)\big),
\,\,\mathrm{for}\,\, k=0,1, \dots,n,
\end{align*}
\end{theorem}

Letting $u=-1$, we have the next theorem as a corollary to Theorem 6.1.

\begin{theorem}
Let $p(x) \in \mathbb{C}[x], with \,\,\mathrm{deg}\, p(x)=n$, and let $f(t)=\frac{1}{\lambda}(e^{\lambda t}-1)$. Then we have 
$p(x)=\sum_{k=0}^{n}a_k \mathcal{E}_{k,\lambda}(x)$, where
\begin{align*}
&a_{k}=\frac{1}{2k!}(f(t))^{k}(p(x+1)+p(x))|_{x=0}\\
&=\frac{1}{2k! \lambda^{k}}\big \langle \big(e^{\lambda t}-1\big)^{k} \big | p(x+1)+p(x)\big \rangle \\
&=\frac{1}{2k! \lambda^{k}}\Delta_{\lambda}^{k}\tilde{\Delta} p(0) \\
&=\frac{1}{2k! \lambda^{k}}\sum_{j=0}^{k}\binom{k}{j}(-1)^{k-j}(p(1+j \lambda)+p(j \lambda)) \\
&=\frac{1}{2}\sum_{l=r}^{n}S_{2}(l,r)\frac{\lambda^{l-r}}{l!}\big(p^{(l)}(1)+p^{(l)}(0)\big),
\,\,\mathrm{for}\,\, k=0,1, \dots,n,
\end{align*}
\end{theorem}

\begin{remark}
Let $p(x) \in \mathbb{C}[x], with \,\,\mathrm{deg}\, p(x)=n$. \par
(a) Write $p(x)=\sum_{k=0}^{n} a_{k}H_{k}(x|u)$. As $\lambda$ tends to $0$, $f(t) \rightarrow t$. Thus we recover from Theorem 3.1 the result in \textnormal{[9,10]}. Namely, we have
\begin{equation}\label{13C}
a_{k}=\frac{1}{(1-u)k!}(p^{(k)}(1)-u p^{(k)}(0)),\,\,\textnormal{for}\,\,k=0,1,\dots,n.
\end{equation}
(b) Write $p(x)=\sum_{k=0}^{n} b_{k}E_{k}(x)$. By letting $u=-1$, we recover from \eqref{13C} the result in \textnormal{[12]}. Namely, we have
\begin{equation*}
b_{k}=\frac{1}{2k!}(p^{(k)}(1)+p^{(k)}(0)),\,\,\textnormal{for}\,\,k=0,1,\dots,n.
\end{equation*}
\end{remark}

\section{Representation by higher-order degenerate Frobenius-Euler polynomials}

Our goal here is to deduce formulas expressing any polynomial in terms of the higher-order degenerate Frobenius-Euler polynomials $h_{n,\lambda}^{(r)}(x|u)$. \par
With $g(t)=\frac{e^t -u}{1-u}, f(t)=\frac{1}{\lambda}(e^{\lambda t}-1)$, from \eqref{10A} and \eqref{11B} we note that
\begin{equation*}
h_{n,\lambda}^{(r)}(x|u) \sim (g(t)^{r}, f(t)).
\end{equation*}
From \eqref{12B}, we have
\begin{equation}\label{1D}
f(t)h_{n,\lambda}^{(r)}(x|u)=nh_{n-1,\lambda}^{(r)}(x|u),
\end{equation}
and from the generating function of the higher-order degenerate Bernoulli polynomials, it is immediate to see that
\begin{equation}\label{2D}
h_{n,\lambda}^{(r)}(x+1|u)-uh_{n,\lambda}^{(r)}(x|u)=(1-u)h_{n,\lambda}^{(r-1)}(x|u).
\end{equation}
As \eqref{2D} is equivalent to saying that $g(t)h_{n,\lambda}^{(r)}(x|u)=h_{n,\lambda}^{(r-1)}(x|u)$,
we have
\begin{equation}\label{3D}
g(t)^{r}h_{n,\lambda}^{(r)}(x|u)=h_{n,\lambda}^{(0)}(x|u)=(x)_{n,\lambda}.
\end{equation}
Now, we assume that $p(x) \in \mathbb{C}[x]$ has degree $n$, and write $p(x)=\sum_{k=0}^{n}a_{k}h_{k,\lambda}^{(r)}(x|u)$.
Then we observe that
\begin{equation}\label{4D}
g(t)^{r}p(x)=\sum_{l=0}^{n}a_{l}\,g(t)^{r}h_{l,\lambda}^{(r)}(x|u)=\sum_{l=0}^{n}a_{l}(x)_{l,\lambda}.
\end{equation}
By using \eqref{4D} and \eqref{4C}, we observe that
\begin{align}
f(t)^{k}g(t)^{r}p(x)&=\sum_{l=0}^{n}a_{l}f(t)^{k}(x)_{l,\lambda}\label{5D}\\
&=\sum_{l=k}^{n}a_{l}\,l(l-1) \cdots (l-k+1)(x)_{l-k,\lambda}.\nonumber
\end{align}
By evaluating \eqref{5D} at $x=0$, we obtain
\begin{equation}\label{6D}
a_k=\frac{1}{k!}f(t)^{k}g(t)^{r}p(x)|_{x=0}= \frac{1}{k!} \langle f(t)^{k}g(t)^{r} |p(x) \rangle .
\end{equation}
This also follows from the observation $\langle g(t)^{r}f(t)^{k} | h_{l,\lambda}^{(r)}(x|u) \rangle=l!\,\delta_{l,k}.$ \par
The equations \eqref{13A} and \eqref{7B} give some alternative expressions of \eqref{6D} as in the following:
\begin{align}
a_k&=\frac{1}{k!}g(t)^{r}f(t)^{k}p(x)|_{x=0}\nonumber \\
&=\frac{1}{(1-u)^{r}k!\lambda^{k}}\sum_{j=0}^{r}\binom{r}{j}(-u)^{r-j}\Delta_{\lambda}^{k}p(j) \label{7D} \\
&=\frac{1}{(1-u)^{r}k!\lambda^{k}}\sum_{j=0}^{r}\sum_{l=0}^{k}\binom{r}{j}\binom{k}{l}(-u)^{r-j}(-1)^{k-l}p(l\lambda+j) \nonumber \\
&=\frac{1}{(1-u)^{r}k!\lambda^{k}}\sum_{j=0}^{r-1}\binom{r-1}{j}(-u)^{r-j-1}\Delta_{\lambda}^{k}(p(x+j+1)-up(x+j))|_{x=0}.\nonumber
\end{align}
Another expression for $a_{k}$ follows from \eqref{16A}:
\begin{align}
a_k&=\frac{1}{(1-u)^{r}k!\lambda^{k}}\sum_{j=0}^{r}\binom{r}{j}(-u)^{r-j}\Delta_{\lambda}^{k}p(x+j)|_{x=0}\nonumber\\
&=\frac{1}{(1-u)^{r}\lambda^{k}}\sum_{j=0}^{r}\binom{r}{j}(-u)^{r-j}\frac{1}{k!}(e^{\lambda t}-1)^{k}p(x+j)|_{x=0}\label{8D}\\
&=\frac{1}{(1-u)^{r}\lambda^{k}}\sum_{j=0}^{r}\sum_{l=k}^{n}\binom{r}{j}(-u)^{r-j}\frac{\lambda^{l}}{l!}S_{2}(l,k)p^{(l)}(j).\nonumber
\end{align}

Summarizing the results so far, from \eqref{7D} and \eqref{8D} we obtain the following theorem.
\begin{theorem}
Let $p(x) \in \mathbb{C}[x], with\,\,\mathrm{deg}\, p(x)=n$. Let $g(t)=\frac{e^{t}-u}{1-u}, f(t)=\frac{1}{\lambda}(e^{\lambda t}-1)$. Then we have 
\begin{equation*}
p(x)=\sum_{k=0}^{n}a_{k}h_{k,\lambda}^{(r)}(x),
\end{equation*}
\end{theorem}
where 
\begin{align*}
a_k&=\frac{1}{k!}g(t)^{r}f(t)^{k}p(x)|_{x=0}\\
&=\frac{1}{(1-u)^{r}k!\lambda^{k}}\sum_{j=0}^{r}\binom{r}{j}(-u)^{r-j}\Delta_{\lambda}^{k}p(j)\\
&=\frac{1}{(1-u)^{r}k!\lambda^{k}}\sum_{j=0}^{r}\sum_{l=0}^{k}\binom{r}{j}\binom{k}{l}(-u)^{r-j}(-1)^{k-l}p(l\lambda+j)\\
&=\frac{1}{(1-u)^{r}k!\lambda^{k}}\sum_{j=0}^{r-1}\binom{r-1}{j}(-u)^{r-j-1}\Delta_{\lambda}^{k}(p(x+j+1)-up(x+j))|_{x=0}\\
&=\frac{1}{(1-u)^{r}\lambda^{k}}\sum_{j=0}^{r}\sum_{l=k}^{n}\binom{r}{j}(-u)^{r-j}\frac{\lambda^{l}}{l!}S_{2}(l,k)p^{(l)}(j).
\end{align*}

We get the next result as a corollary to Theorem 4.1 by letting $u=-1$.

\begin{theorem}
Let $p(x) \in \mathbb{C}[x], \mathrm{deg}\, p(x)=n$. Let $g(t)=\frac{e^{t}+1}{2}, f(t)=\frac{1}{\lambda}(e^{\lambda t}-1)$. Then we have 
\begin{equation*}
p(x)=\sum_{k=0}^{n}a_{k}\mathcal{E}_{k,\lambda}^{(r)}(x),
\end{equation*}
\end{theorem}
where 
\begin{align*}
a_k&=\frac{1}{k!}g(t)^{r}f(t)^{k}p(x)|_{x=0}\\
&=\frac{1}{2^{r}k!\lambda^{k}}\tilde{\Delta}^{r}\Delta_{\lambda}^{k}p(0)\\
&=\frac{1}{2^{r}k!\lambda^{k}}\sum_{j=0}^{r}\sum_{l=0}^{k}\binom{r}{j}\binom{k}{l}(-1)^{k-l}p(j+l \lambda)\\
&=\frac{1}{2^{r}k!\lambda^{k}}\tilde{\Delta}^{r-1}\Delta_{\lambda}^{k}(p(x+1)+p(x))|_{x=0} \\
&=\frac{1}{2^{r}\lambda^{k}}\sum_{j=0}^{r}\sum_{l=k}^{n}\binom{r}{j}\frac{\lambda^{l}}{l!}S_{2}(l,k)p^{(l)}(j).
\end{align*}

\begin{remark}
Let $p(x) \in \mathbb{C}[x], with \,\,\mathrm{deg}\, p(x)=n$. \par
(a) Write $p(x)=\sum_{k=0}^{n}a_k H_{k}^{(r)}(x|u)$. As $\lambda$ tends to $0$, $f(t) \rightarrow t$. Thus, from Theorem 4.1, we recover the following result obtained in \textnormal{[9,10]}: \par
\begin{align}
a_k&=\frac{1}{k!}g(t)^{r}f(t)^{k}p(x)|_{x=0}\nonumber\\
&=\frac{1}{(1-u)^{r}k!}\sum_{j=0}^{r}\binom{r}{j}(-u)^{r-j}p^{(k)}(j)\label{9D}
\end{align}
\begin{align*}
&=\frac{1}{(1-u)^{r}k!}\sum_{j=0}^{r-1}\binom{r-1}{j}(-u)^{r-j-1}(p^{(k)}(j+1)-up^{(k)}(j)),
\end{align*}
where $g(t)=\frac{e^{t}-u}{1-u},\,and \,\, f(t)=t$. \par
(b) Write $p(x)=\sum_{k=0}^{n}b_k E_{k}^{(r)}(x)$. Now, by letting $u=-1$ in \eqref{9D}, we recover the following result obtained in \textnormal{[8]}: \par
\begin{align*}
b_k&=\frac{1}{k!}g(t)^{r}f(t)^{k}p(x)|_{x=0}\\
&=\frac{1}{2^{r}k!}\sum_{j=0}^{r}\binom{r}{j}p^{(k)}(j)\\
&=\frac{1}{2^{r}k!}\sum_{j=0}^{r-1}\binom{r-1}{j}(p^{(k)}(j+1)+p^{(k)}(j)),
\end{align*}
where $g(t)=\frac{e^{t}+1}{2},\,and \,\, f(t)=t$.
\end{remark}

\section{Examples}
(a) Here we illustrate Theorem 3.1, with $p(x)=H_{n}(x|u)$. Let $H_{n}(x|u)=\sum_{k=0}^{n}a_{k}h_{k,\lambda}(x|u)$.\par
For $k$ with $0  \le k \le n$, we have
\begin{align}
a_{k}&=\frac{1}{(1-u)k! \lambda^{k}}\big \langle \big(e^{\lambda t}-1\big)^{k} \big | H_{n}(x+1|u)-uH_{n}(x|u)\big \rangle \nonumber\\
&=\frac{1}{k! \lambda^{k}}\big \langle \big(e^{\lambda t}-1\big)^{k} \big | x^{n}\big \rangle \label{1E}\\
&=\frac{1}{k! \lambda^{k}}\Delta_{\lambda}^{k}x^{n}|_{x=0}\nonumber\\
&=\frac{1}{k! \lambda^{k}}\Delta_{\lambda}^{k}0^{n}.\nonumber
\end{align}
Noting that $H_n(x|u) \sim \big(\frac{e^{t}-u}{1-u},t \big),\,\, h_{n,\lambda}(x|u) \sim \big(\frac{e^t -u}{1-u}, \frac{1}{\lambda}(e^{\lambda t}-1)\big)$, one can determine the coefficients also by using \eqref{15B} and \eqref{16B}. Thus, from \eqref{1E} we have 
\begin{equation*}
H_{n}(x|u)=\sum_{k=0}^{n}\frac{1}{k! \lambda^{k}}\Delta_{\lambda}^{k}0^{n}h_{k,\lambda}(x|u).
\end{equation*}
(b) Here we illustrate Theorem 3.1, for $p(x)=\sum_{k=1}^{n-1}\frac{1}{k(n-k)}B_{k}(x)B_{n-k}(x),\,\,(n \ge 2)$. For this, we first recall from [15] that
\begin{align*}
p(x)&=\sum_{k=0}^{n-2}
\bigg(\frac{2\binom{n}{k}}{n} \bigl(H_{n-1} - H_{n-k-1}\bigl)B_{n-k} + \binom{n-1}{k}  \sum_{l=k+1}^{n-1} \frac{B_{l-k}B_{n-l}}{(l-k)(n-l)}
\\&\quad \quad\quad+ \frac{\binom{n-1}{k} B_{n-1-k}}{n-1-k} \bigg)E_{k}(x) +\frac{n-1}{2} E_{n-2}(x) + \frac{2}{n}H_{n-1}E_{n}(x),
\end{align*}
where $H_{n}=1+\frac{1}{2}+\cdots+\frac{1}{n}$ is the harmonic number. 
Let $p(x)=\sum_{k=0}^{n}a_{k}
\mathcal{E}_{k,\lambda}(x)$. 
For $k$, with $ 0 \le k \le n$, we obtain
\begin{align}
2k! \lambda^{k}a_{k}&=\big \langle \big(e^{\lambda t}-1\big)^{k} \big | p(x+1)+p(x)\big \rangle \nonumber\\
&=\sum_{m=0}^{n-2}
\bigg(\frac{2\binom{n}{m}}{n} \bigl(H_{n-1} - H_{n-m-1}\bigl)B_{n-m} + \binom{n-1}{m}  \sum_{l=m+1}^{n-1} \frac{B_{l-m}B_{n-l}}{(l-m)(n-l)}\nonumber\\
&\quad \quad\quad+ \frac{\binom{n-1}{m} B_{n-1-m}}{n-1-m} \bigg)\big \langle \big(e^{\lambda t}-1\big)^{k} \big | E_{m}(x+1)+E_{m}(x)\big \rangle
\label{2E}\\
&\quad \quad\quad+\frac{n-1}{2}\langle \big(e^{\lambda t}-1\big)^{k} \big | E_{n-2}(x+1)+E_{n-2}(x)\big \rangle\nonumber\\
&\quad\quad\quad+ \frac{2}{n}H_{n-1}\langle \big(e^{\lambda t}-1\big)^{k} \big | E_{n}(x+1)+E_{n}(x)\big \rangle. \nonumber
\end{align}
By proceeding as we did in (a), from \eqref{2E} we see that
\begin{align}
k! \lambda^{k}a_{k}&=\sum_{m=0}^{n-2}
\bigg(\frac{2\binom{n}{m}}{n} \bigl(H_{n-1} - H_{n-m-1}\bigl)B_{n-m} + \binom{n-1}{m}  \sum_{l=m+1}^{n-1} \frac{B_{l-m}B_{n-l}}{(l-m)(n-l)}\label{3E}\\
&\quad \quad\quad+ \frac{\binom{n-1}{m} B_{n-1-m}}{n-1-m} \bigg) \Delta_{\lambda}^{k}0^{m}
+\frac{n-1}{2} \Delta_{\lambda}^{k}0^{n-2}+ \frac{2}{n}H_{n-1}\Delta_{\lambda}^{k}0^{n}. \nonumber
\end{align}

Thus, from \eqref{3E} and for $ n \ge 2$, we have
\begin{align*}
&\sum_{k=1}^{n-1}\frac{1}{k(n-k)}B_{k}(x)B_{n-k}(x)\\
&=\sum_{k=0}^{n}\frac{1}{k! \lambda^{k}}\bigg\{\sum_{m=0}^{n-2}
\bigg(\frac{2\binom{n}{m}}{n} \bigl(H_{n-1} - H_{n-m-1}\bigl)B_{n-m} + \binom{n-1}{m}  \sum_{l=m+1}^{n-1} \frac{B_{l-m}B_{n-l}}{(l-m)(n-l)}\nonumber\\
&\quad \quad\quad+ \frac{\binom{n-1}{m} B_{n-1-m}}{n-1-m} \bigg) \Delta_{\lambda}^{k}0^{m}
+\frac{n-1}{2} \Delta_{\lambda}^{k}0^{n-2}+ \frac{2}{n}H_{n-1}\Delta_{\lambda}^{k}0^{n}\bigg\}\mathcal{E}_{k,\lambda}(x).
\end{align*}

(c) In [15], it is shown that the following identity holds for $n \ge 2$:
\begin{align*}
\sum_{k=1}^{n-1}&\frac{1}{k(n-k)}E_{k}(x)E_{n-k}(x)\\
&=\sum_{k=0}^{n-2}\sum_{l=k+1}^{n-1}
\frac{\binom{n-1}{k}}{(l-k)(n-l)}E_{l-k}E_{n-l}E_{k}(x)+\frac{2}{n}H_{n-1}E_{n}(x),
\end{align*}
where $H_{n}=1+\frac{1}{2}+\cdots+\frac{1}{n}$ is the harmonic number. \par
Write $\sum_{k=1}^{n-1} \frac{1}{k(n-k)}E_{k}(x) E_{n-k} (x)=\sum_{k=0}^{n}a_{k}\mathcal{E}_{k,\lambda}(x)$. \par
By proceeding similarly to (b), we obtain the following identity:
\begin{align*}
&\sum_{k=1}^{n-1} \frac{1}{k(n-k)}E_{k}(x) E_{n-k} (x) \\
&=\sum_{k=0}^{n}\frac{1}{k! \lambda^{k}}\bigg(\sum_{m=0}^{n-2}\sum_{l=m+1}^{n-1}
\frac{\binom{n-1}{m}}{(l-m)(n-l)}E_{l-m}E_{n-l}\Delta_{\lambda}^{k}0^{m}+\frac{2}{n}H_{n-1}\Delta_{\lambda}^{k}0^{n}\bigg)\mathcal{E}_{k,\lambda}(x).
\end{align*}

(d) In [15], it is proved that the following identity is valid for $n \ge 2$:
\begin{align*}
&\sum_{k=1}^{n-1}\frac{1}{k(n-k)}B_{k}(x)E_{n-k}(x)
\\
&=\sum_{k=0}^{n-2}\Bigg(\frac{1}{n}\binom{n}{k}  \big (H_{n-1} -H_{n-k-1} \big) B_{n-k} - \frac{1}{2}\binom{n-1}{k} \frac{E_{n-k-1}}{n-k-1}\Bigg)E_{k}(x)
+\frac{2}{n}H_{n-1}E_{n}(x).
\end{align*}
Again, by proceeding analogously to (b), we get the following identity:
\begin{align*}
&\sum_{k=1}^{n-1}\frac{1}{k(n-k)}B_k(x)E_{n-k}(x) \\
&=\sum_{k=0}^{n}\frac{1}{k! \lambda^{k}}\bigg\{\sum_{m=0}^{n-2}\Bigg(\frac{1}{n}\binom{n}{m}  \big (H_{n-1} -H_{n-m-1} \big) B_{n-m} - \frac{1}{2}\binom{n-1}{m} \frac{E_{n-m-1}}{n-m-1}\Bigg)\Delta_{\lambda}^{k}0^{m}\\
&\quad\quad\quad+\frac{2}{n}H_{n-1}\Delta_{\lambda}^{k}0^{n}\bigg\}\mathcal{E}_{k,\lambda}(x).
\end{align*}

(e) Nielsen [18,2] expressed the product of a Bernoulli polynomial and a Frobenius-Euler polynomial in terms of Frobenius-Euler polynomials. Namely, for nonnegative integers $m$ and $n$, 
\begin{align*}
&B_{m}(x)H_{n}(x|u)\\
&=mH_{m+n-1}(x|u)+\sum_{r=0}^{m}\binom{m}{r}B_{r}H_{m+n-r}(x|u)+\frac{mu}{1-u}\sum_{s=0}^{n}\binom{n}{s}H_{s}(u)H_{m+n-s-1}(x|u).
\end{align*}

Again, in a similar way to (b), we can show that
\begin{align*}
&B_{m}(x)E_{n}(x) \\
&=\sum_{k=0}^{m+n}\frac{1}{k! \lambda^{k}}\bigg\{m\Delta_{\lambda}^{k}0^{m+n-1}+\sum_{r=0}^{m}\binom{m}{r}B_{r}\Delta_{\lambda}^{k}0^{m+n-r}\\
&\quad\quad\quad\quad\quad\quad+\frac{mu}{1-u}\sum_{s=0}^{n}\binom{n}{s}H_{s}(u)\Delta_{\lambda}^{k}0^{m+n-s-1}\bigg\}h_{k,\lambda}(x).
\end{align*}
We also represent $B_{m}(x)H_{n}(x|u)$ in terms of the degenerate Frobenius-Euler polynomials of order $r$. Indeed, from Theorem 4.1 we obtain two expressions:\\
\begin{align*}
B_{m}(x)H_{n}(x|u)&=\sum_{k=0}^{m+n}\frac{1}{(1-u)^{r}k!\lambda^{k}}\bigg(\sum_{j=0}^{r}\sum_{l=0}^{k}\binom{r}{j}\binom{k}{l}(-u)^{r-j}(-1)^{k-l} \\
&\quad\quad\quad \times B_{m}(j+l \lambda)H_{n}(j+l \lambda|u)\bigg)h_{k,\lambda}^{(r)}(x|u) \\
&=\sum_{k=0}^{m+n}\frac{1}{(1-u)^{r-1}k!\lambda^{k}}\sum_{j=0}^{r-1}\binom{r-1}{j}(-u)^{r-j-1}\bigg\{m\Delta_{\lambda}^{k}0^{m+n-1}\\
&+\sum_{r=0}^{m}\binom{m}{r}B_{r}\Delta_{\lambda}^{k}0^{m+n-r}
+\frac{mu}{1-u}\sum_{s=0}^{n}\binom{n}{s}H_{s}(u)\Delta_{\lambda}^{k}0^{m+n-s-1}\bigg\}h_{k,\lambda}^{(r)}(x|u).
\end{align*}

(f) In [2], Carlitz showed the following identity: for $u \ne 1, v \ne 1, and\,\, uv \ne 1$, we have
\begin{align*}
H_{m}(x|u)H_{n}(x|v)&=H_{m+n}(x|uv) \\
&+\frac{u(1-v)}{1-uv}\sum_{r=1}^{m}\binom{m}{r}H_{r}(u)H_{m+n-r}(x|uv) \\
&+\frac{v(1-u)}{1-uv}\sum_{s=1}^{n}\binom{n}{s}H_{s}(v)H_{m+n-s}(x|uv).
\end{align*}
Write $H_{m}(x|u)H_{n}(x|v)=\sum_{k=0}^{m+n} a_{k}h_{k,\lambda}(x|uv)$. Then, proceeding similarly to (b), we have:
\begin{align*}
H_{m}(x|u)H_{n}(x|v)&=\sum_{k=0}^{m+n}\frac{1}{k!\lambda^{k}}\bigg\{\Delta_{\lambda}^{k}0^{m+n}+\frac{u(1-v)}{1-uv}\sum_{r=1}^{m}\binom{m}{r}H_{r}(u)\Delta_{\lambda}^{k}0^{m+n-r}\\
&+\frac{v(1-u)}{1-uv}\sum_{s=1}^{n}\binom{n}{s}H_{s}(v)\Delta_{\lambda}^{k}0^{m+n-s}\bigg\}h_{k,\lambda}(x|uv).
\end{align*}

\section{Conclusion}

In this paper, we were interested in representing any polynomial in terms of the degenerate Frobenius-Euler polynomials and of the higher-order degenerate Frobenius-Euler polynomials. We were able to derive formulas for such representations with the help of umbral calculus. We showed that,  by letting $\lambda$ tend to zero, they agree with the previously found formulas for representations by the Frobenius-Euler polynomials and by the higher-order Frobenius-Euler polynomials. In addition, by letting $u=-1$ we obtained formulas for expressing any polynomial by the degenerate Euler polynomials and by the higher-order degenerate Euler polynomials. Further, we illustrated the formulas with some examples. \\
Even though the method adopted in this paper is simple, they are very useful and powerful. Indeed, as we mentioned in the Section 1, both a variant of Miki's identity and Faber-Pandharipande-Zagier (FPZ) identity follow from the one identity (see \eqref{1A}) that can be derived from a formula (see [12]) involving only derivatives and integrals of the given polynomial, while all the other known proofs are quite involved. \\
It is one of our future research projects to continue to find formulas for representing polynomials in terms of some specific special polynomials and to apply those to discovering some interesting identities.


\end{document}